\documentclass[11pt,a4paper]{article}
\usepackage[a4paper,top=3cm,bottom=4.5cm,left=2.5cm,right=2.5cm]{geometry}
\usepackage{graphicx,caption}
\usepackage{amssymb}
\usepackage{mathrsfs}
\usepackage{amsmath}
\usepackage{latexsym}
\usepackage{amsthm}
\usepackage{eucal}
\usepackage{eufrak}
\usepackage{relsize}
\numberwithin{equation}{section}
\newtheorem{teo}{Theorem}
\newtheorem{defn}{Definition}

\theoremstyle{definition}

\begin{document}
	\begin{center}
		\Large\textbf{Torsional instability and sensitivity analysis in a suspension bridge model related to the Melan equation}\\
		\vspace{4mm}
		\large{Alessio Falocchi}
		\vspace{5mm}
		\begin{center}
			{\scriptsize Dipartimento di Matematica - Politecnico di Milano\\	Piazza Leonardo da Vinci 32 - 20133 Milano, Italy\\	\texttt{alessio.falocchi@polimi.it}}
		\end{center}
	\end{center}
	\vspace{3mm}
	\begin{abstract}
		Inspired by the Melan equation we propose a model for suspension bridges with two cables linked to a deck, through inextensible hangers. We write the energy of the system and we derive from variational principles two nonlinear and nonlocal hyperbolic partial differential equations, involving the vertical displacement and the torsional rotation of the deck. We prove existence and uniqueness of a weak solution and we perform some numerical experiments on the isolated system; moreover we propose a sensitivity analysis of the system by mechanical parameters in terms of torsional instability. Our results display that there are specific thresholds of torsional instability with respect to the initial amplitude of the longitudinal mode excited.
	\end{abstract}
	
	\section{Introduction}
	\label{Introduction}
	The Melan equation was introduced by the Austrian engineer Josef Melan \cite{melan} in 1888 to model a suspension bridge; Melan considered the bridge as a combination of a string (the cable) and a beam (the deck) linked through some rigid hangers, which are considered uniformly distributed along the main span.\\  The equation can be derived writing the equilibrium of the beam and  the  string and combining the two equations through the live load, carried in part by the cable and in part by the deck. The result is the following fourth order differential equation   
	\begin{equation}
	\begin{cases}
	EIw''''(x)-(H+h(w))w''(x)-\frac{q}{H}h(w)=p(x) & \forall x\in (0,L)\\
	w(0)=w(L)=w''(0)=w''(L)=0,\\
	\end{cases}
	\label{Melan}
	\end{equation}
	in which $w(x)$ is the vertical displacement of the beam (positive if directed downward), $EI$ is the flexural rigidity of the beam, $H$ is the horizontal tension of the string when subjected to the dead load $-q$, and $h(w)$ is a nonlocal term, representing the additional tension in the cable due to the live load $p(x)$; the beam has a span equal to $L$ and is supposed hinged at the endpoints. 
	\\ The presence of the nonlocal term makes challenging the study of the equation from both the theoretical as from the numerical point of view, see e.g. \cite{samet,wang,semper}; although  (\ref{Melan}) cannot be derived from the variation of the corresponding energy \cite{wang}, von K\'{a}rm\'{a}n-Biot \cite{biot} call the Melan equation (\ref{Melan}) the \textit{fundamental equation of the theory of the suspension bridge.} \\
	This equation is our starting point, we propose a more reliable model for suspension bridge in which we have two strings (the cables) linked to the same deck, through inextensible hangers, see Section \ref{Description of the dynamical model}. In this way we introduce the torsional rotation of the deck, which cannot be seen in a one-dimensional model. Our two main purposes are to study the torsional instability and to analyze how the mechanical parameters of the bridge affect this behavior; these purposes are motivated by the fact that in suspension bridges torsional oscillations can be catastrophic, see for instance the video of the collapse of the Tacoma Narrows bridge (TNB) \cite{Video}. This case is not isolated, in \cite[pp.1-40]{Gazz} the author mentions many other suspension bridges that manifested this behavior. For instance, the Brighton Chain Pier collapsed in 1833 due to windstorms that caused different kinds of oscillation included the torsional ones; similar circumstances occurred for the Menai Straits Bridge in 1839 and the Wheeling Suspension Bridge in West Virginia (1854), where the witnesses told about a "twisted and writhed" movement that lasted only "two minutes". \\After the TNB collapse many hypotheses have been proposed to explain the so-called torsional instability; some explanations were found in the aerodynamic effects, from the vortex shedding to the parametric resonance and the flutter theory. Nowadays there are still many doubts and new suspension bridges continue to manifest dangerous and sudden oscillations. Matukituki Suspension Footbridge collapsed in 1977, showing a noticeable node at midspan typical of the torsional motion; in 2000 all the world spoke about the closure of the London Millenium Bridge two days after its inauguration, since the crowd passing over it caused strange vibrations. For further details and other events of such type we refer to \cite{Gazz}; it is interesting to see that in the most of the cases there are suspension bridges that suddenly change harmless vertical oscillations, possibly due to the wind or the pedestrians' walk, into different dangerous oscillations. In this paper we propose a model for suspension bridges able to catch this activation phenomenon thanks to the nonlinear configuration of the structure.\\  
	In Section \ref{Energies involved in the structure} we compute the energies involved in the system and we derive the corresponding Euler-Lagrange equations by variational principles. We obtain a system of two nonlinear partial differential equations in space and time with nonlocal terms, as in the original Melan equation, see Section \ref{The Euler-Lagrange Equations}. The nonlinearities are due to the geometric configuration of the suspension bridge, that involves the parabolic shape of the cables, see also \cite{lacarbonara}, and the rotation of the deck, in which trigonometric functions appear; the linearization of trigonometric
	functions is admissible assuming small torsional rotation, but we avoid it, complying with the real bridge geometry and the possible presence of large rotations, about this topic see also \cite{falocchi}.\\
	The existence and uniqueness of a weak solution in the proper functional spaces is proved in Section \ref{Proof of the Theorem}, applying the Galerkin procedure to the hyperbolic equations; we give the complete version of the proof because the presence of the nonlinearities makes challenging the uniqueness problem, that is proved in a wider functional space. In particular we obtain the latter result in a non-standard way, testing the equations with the "potential", i.e. with the Green function applied to the time derivative of the solutions.\\We consider an isolated model aiming to show that the origin of the torsional
	instability is purely structural, as proposed in other works \cite{Gazzola hyperb, Gazzola Torsion,falocchi}; in particular we suppose that the wind introduces energy in the structure, exciting one longitudinal mode at a time, through the initial conditions. This is legitimate since the frequency of the vortex shedding usually excites one mode. Our numerical results, presented in Section \ref{Numerical experiments}, show that there exist thresholds of torsional instability, depending on the initial amplitude of the longitudinal mode excited. \\A discussion about the influence of the mechanical parameters of the bridge on the torsional stability is provided in Section \ref{The influence of the mechanical parameters on the stability of the system}, giving some hints to bridge designers; in the analysis of sensitivity we show that an important role on the stability is assumed by the sag-span ratio, the shape of the cross section of the deck and its mass.\\  
	\section{The suspension bridge model}
	\label{The dynamical model}
	\subsection{Description of the dynamical model}
	\label{Description of the dynamical model}
	The main elements composing a suspension bridge are four towers, a rectangular deck, two sustaining cables and some hangers; the cables, fixed to the towers, sustain the deck through inextensible hangers, as in \cite{luco}. 
	\begin{figure}[h]
		\centering
		\includegraphics[width=11cm]{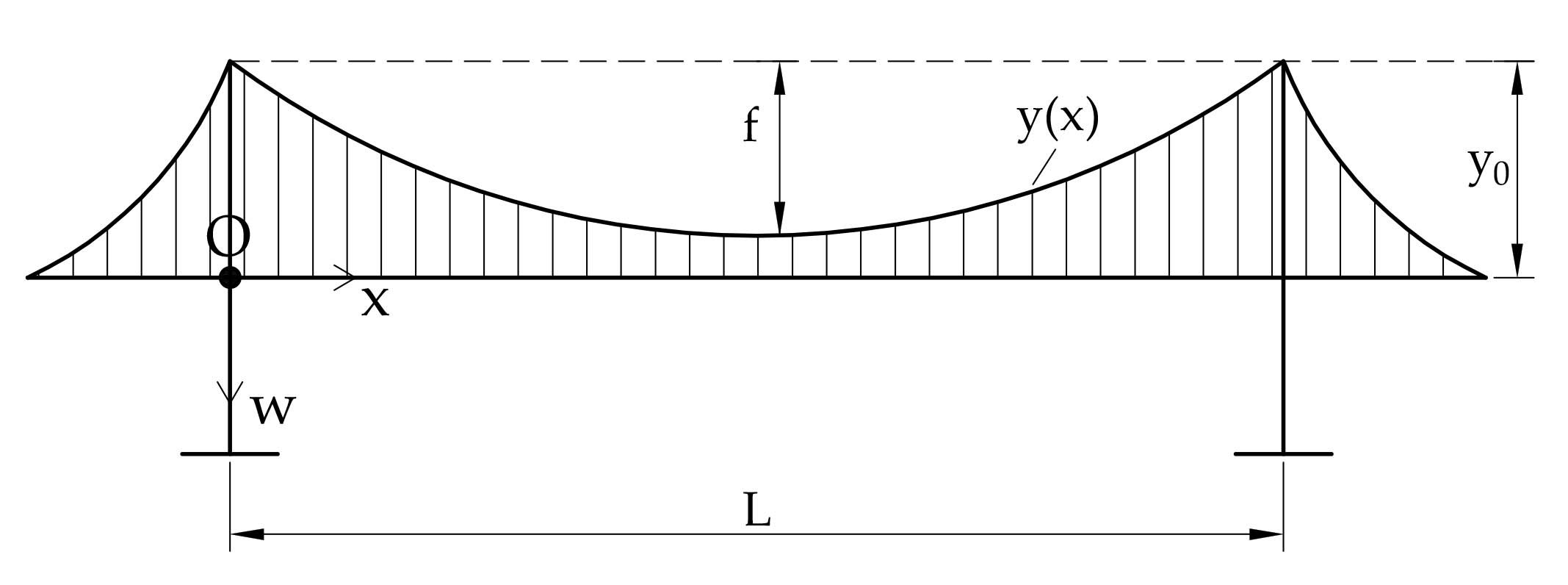}
		\caption{Sketch of the side view of the suspension bridge with the quotes assumed positive.}
		\label{ponte}
	\end{figure}
	The deck is modeled as a degenerate plate, composed by a beam with length $L$, corresponding to its midline, and cross sections with length $2\ell\ll L$; the midline connects the barycentres of the cross sections and the latter may rotate with respect to the horizontal equilibrium position. The edges of the degenerate plate, i.e. the endpoints of the cross sections, are linked to the cables by the hangers, see also \cite{berchio,falocchi}.\\
	We introduce a reference system $(O,x,y)$ with origin in correspondence of a tower at the level of the deck, assuming $w$ positive if directed downward and $x$ along the main span of the bridge, see Figure \ref{ponte}. \\
	Following Melan \cite{melan} and von K\'{a}rm\'{a}n-Biot \cite[Section VII.5]{biot}, we consider the main span of the suspension bridge as a combined system of two perfectly flexible strings (the main cables) linked to the deck through inextensible hangers; this model is more complex than the original one of von K\'{a}rm\'{a}n-Biot, because in the latter there was only a string linked to a single beam. In our model, instead, a movement of the deck influences both the cables and the result is a system of equations strongly coupled. \\
	We denote the derivatives of a function $f=f(t)$, depending only on $t$, and of a function $g=g(x)$, depending only on $x$, respectively by
	\begin{equation*}
	\dot{f}=\dfrac{df}{dt}, \hspace{6mm} g'=\dfrac{dg}{dx},
	\end{equation*}
	while we denote the partial derivatives of a function $w=w(x,t)$ by
	\begin{equation*}
	w_{x}=\frac{\partial w}{\partial x}, \hspace{4mm} w_{t} =\frac{\partial w}{\partial t}
	\end{equation*} 
	and in a similar way higher order derivatives.\\
	We suppose that the two main cables have the same mechanical properties and that the hangers are uniformly distributed along the two free edges of the deck.
	As suggested in \cite{biot}, we assume that the cables at rest take the shape $y(x)$, where $y$ solves the following differential equation 
	\begin{equation}
	\begin{cases}
	Hy''(x) =-q & \forall x\in (0,L)\\
	y(0)=y(L)=-y_{0}&(y_0>0)\\
	\end{cases}
	\label{eq cavo}
	\end{equation}\\
	Here, $q$ is the dead load, $H$ is the tension of the cable, $L$ is length of the bridge span and $y_{0}$ is the height of the towers. We assume that, for the two cables, the dead load in the initial configuration is $q=\frac{Mg}{2}$, where $M$ is the linear density of the deck mass, $g$ is the gravitational acceleration.
	Since $q$ is constant we find that the cables have the shape of a parabola given by
	\begin{equation}
	y(x)=-\dfrac{Mg}{4H}x^2+\dfrac{MgL}{4H}x-y_0 \hspace{7mm}\forall x\in (0,L)
	\label{y(x)}
	\end{equation}
	and we have 
	\begin{equation*}
	y'(x)=\dfrac{Mg}{2H}\bigg(\dfrac{L}{2}-x\bigg),\hspace{7mm} y''(x)=-\dfrac{Mg}{2H} \hspace{7mm}\forall x\in (0,L).
	\label{y'(x)}
	\end{equation*}
	As suggested in \cite[p.59]{Podolny}, from the elastic theory the parabolic shape of the cables implies that, in a situation of equilibrium, 
	\begin{equation}
	H=\frac{qL^2}{8f}=\frac{MgL^2}{16f},
	\label{H}
	\end{equation}
	where $f$ is the cable sag as in Figure \ref{ponte}. 
	Then an equivalent way to write (\ref{y(x)}) is
	\begin{equation*}
	y(x)=-\dfrac{4f}{L^2}x^2+\dfrac{4f}{L}x-y_0 \hspace{7mm}\forall x\in (0,L).
	\end{equation*}
	Let us introduce the bounded function for all $x\in(0,L)$
	\begin{equation}
	\xi(x):=\sqrt{1+y'(x)^2}, \hspace{15mm}1\leq\xi(x)<\overline{\xi}:=\sqrt{1+\bigg(\dfrac{MgL}{4H}\bigg)^2}
	\label{xi}
	\end{equation}
	which will appear in the calculations.\\
	The deck's deformations in the model comes up, as for a beam, from bending and torsion, due to some energy input; according to the de Saint Venant theory a simple beam has a bending stiffness depending on $E$, the Young modulus, and $I$, the linear density of the moment of inertia of the cross section. On the other hand the beam opposes to torsional movements with a torsional stiffness proportional to $G$, the shear modulus, and $K$, the torsional constant of the section; we point out that the pure torsion, depending on the $GK$-term, occurs only when the warping can take place freely. The presence of welding at the supports,  typical of steel structures, changes in the beam section or imposed torsional moment restrain the warping in some points of the beam \cite{usainstitute}. In 1940, Vlasov \cite{Vlasov} developed a torsional theory in which constrained warping was included; in particular, to the pure torsional term of de Saint Venant, Vlasov added a new differential term of the fourth order, proportional to $E$ and $J$, the warping constant of the section. \\ Since we are considering a model for suspension bridges, whose main structural elements are in steel,  we adopt the Vlasov theory, including in the torsional equation the warping term, to obtain more precise results.\\ About the loading conditions of the bridge we consider the dead load $M$, representing the linear density of the deck, and we neglect the cable and hangers weight, since it is small compared to $M$. 
	\subsection{Energy involved in the structure}
	\label{Energies involved in the structure}
\begin{figure}[h]
		\centering
		\includegraphics[width=10cm]{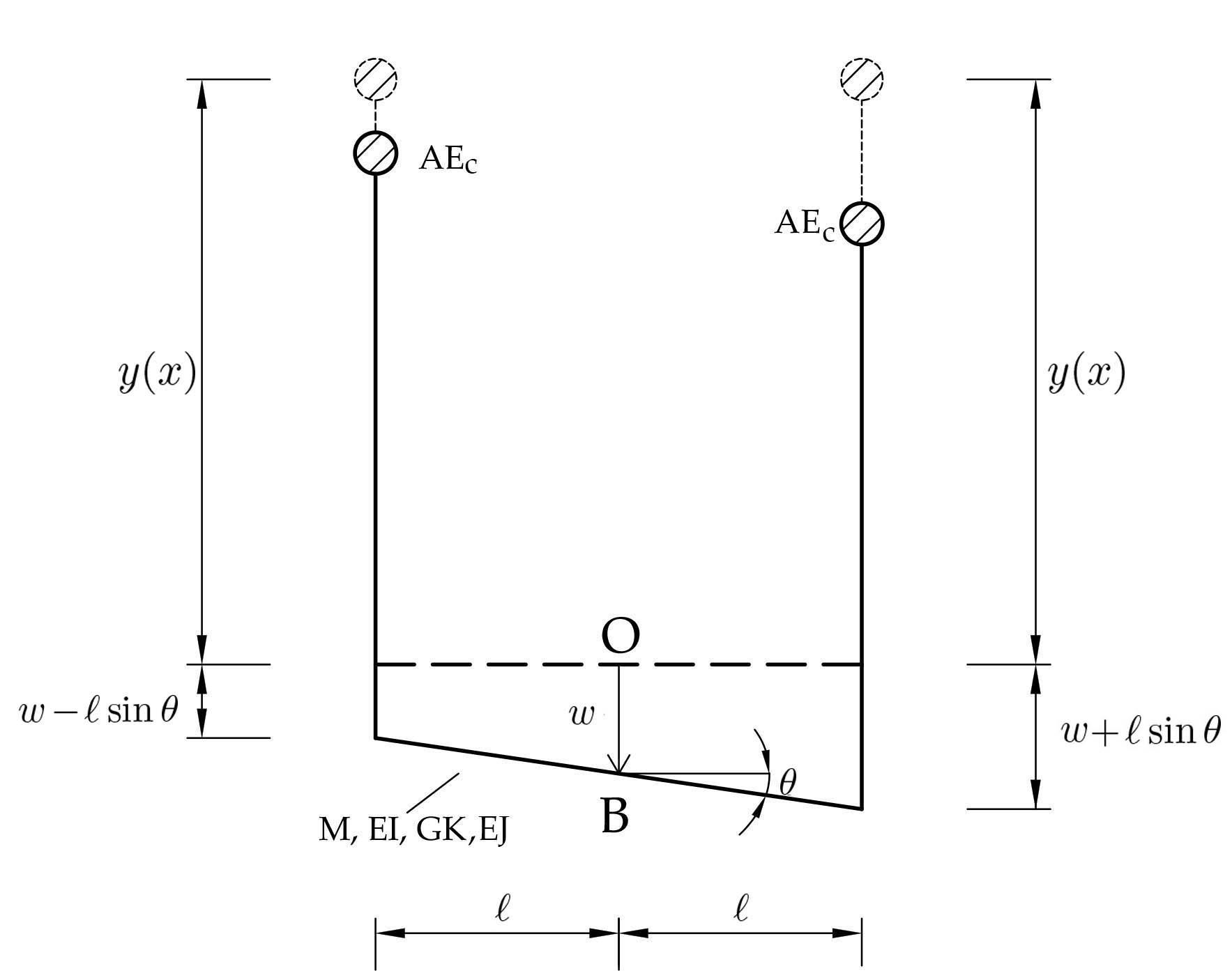}
		\caption{A cross section of the bridge.}
		\label{schema}
	\end{figure}

	In Figure \ref{schema} we sketch a generic cross section of the bridge, highlighting the degrees of freedom of the system given by $w(x,t)$ and $\theta(x,t)$, representing respectively the downward displacement and the torsional rotation of the barycentric line of the deck.\\
	Following the energetic approach suggested in \cite{Gazzola hyperb,falocchi} we compute the energy involved in the cable-hangers-beam system from which we will derive the PDE system. All the mechanical constants are explained in Section \ref{Description of the dynamical model}, then we do not repeat them.
	\begin{itemize}
		\item \textbf{Kinetic energy of the system:} It is obtained adding the vertical and the rotational contribute
		\begin{equation*}
		E_{k}=\frac{M}{2}\int_{0}^{L}w_{t}^2\hspace{1mm}dx+\frac{M\ell^2}{6}\int_{0}^{L}\theta_{t}^2\hspace{1mm}dx.
		\label{en cinetica}
		\end{equation*}
		\item \textbf{Potential energy due to dead loads:} 
		\begin{equation*}
		E_{p}=-Mg\int_{0}^{L} w\hspace{1mm}dx
		\label{en grav deck}
		\end{equation*}
		Since we consider a reference system with $w$ positive downwards this energy is negative. 
		\item \textbf{Stiffening energy of the deck:} It is given by the sum of the bending energy of the deck
		\begin{equation*}
		E_{B}=\frac{EI}{2}\int_{0}^{L}(w_{xx})^2dx
		\label{en bending}
		\end{equation*}
		and the torsional energy 
		\begin{equation*}
		E_{T}=\frac{GK}{2}\int_{0}^{L}(\theta_{x})^2dx+\frac{EJ}{2}\int_{0}^{L}(\theta_{xx})^2dx,
		\label{en tors}
		\end{equation*}
		given by de Saint Venant and Vlasov contributes.
		
		\item \textbf{Stretching energy of the cable:} First of all let us introduce the functional $\Gamma: C^1[0,L]\rightarrow \mathbb{R}$, representing the variation of the length of the cables 
		\begin{equation}
		\begin{split}
		&u \mapsto\Gamma(u):=\int_{0}^{L}\big(\sqrt{1+[(u+y)_{x}]^2}-\sqrt{1+(y')^2}\big)dx\\&\hspace{17.5mm}=\int_{0}^{L}\big(\sqrt{1+[(u+y)_{x}]^2}\big)dx-L_c,
		\end{split}
		\label{gamma}
		\end{equation}
		where $L_c=\int_{0}^{L}\sqrt{1+(y')^2}dx$ is the cable length in the initial configuration.\\
		The tension of the cable is composed by two contributes, the tension at rest
		\begin{equation}
		H(x)=H\xi(x)
		\label{tension}
		\end{equation}
		and the additional tension due to the increment of the length $\Gamma(u)$ of each cable
		\begin{equation}
		\frac{AE_c}{L_{c}}\Gamma(u)
		\label{additional tension}
		\end{equation}
		in which $H$ is the horizontal tension, $A$ the sectional area, $E_c$ the Young modulus of the cable.
		To this terms correspond respectively two deformation energies; the amount of energy needed to deform the cable at rest under the tension (\ref{tension}) from the original position $y(x)$ to $y(x)+u(x,t)$ 
		\begin{equation*}
		E_{C1}(u)=H\int_{0}^{L}\xi(x)\big(\sqrt{1+[(u+y)_{x}]^2}-\xi(x)\big)dx
		\label{energy Hc nonlin}
		\end{equation*}
		and the energy due to the additional tension (\ref{additional tension}) 
		\begin{equation*}
		E_{C2}(u)=	\frac{AE_c}{2L_{c}}\Gamma(u)^2=	\frac{AE_c}{2L_{c}}\bigg(\int_{0}^{L}\big(\sqrt{1+[(u+y)_x]^2}-\xi(x)\big)dx\bigg)^2.
		\label{energy Tc}
		\end{equation*}
		Recalling that the hangers are assumed to be inextensible, from Figure \ref{schema} we see that the vertical displacements of the cables are $u(x,t)=w(x,t)\pm\ell\sin\theta(x,t)$ with respect to the cable considered; then, for a cable, we obtain the following energy
		\begin{equation}
		\begin{split}
		E_{C}(w,\theta)=H\int_{0}^{L}\xi(x)\big(\sqrt{1+[(w+\ell\sin\theta+y)_{x}]^2}-\xi(x)\big)\hspace{1mm}dx+\frac{AE_c}{2L_{c}}[\Gamma(w+\ell\sin\theta)]^2.
		\end{split}
		\label{energy cavi}
		\end{equation}
		By computing the variation of the energy (\ref{energy cavi}) with respect to $w$ and integrating by parts, we find for all $\varphi\in C^{\infty}_c(0,L)$
		\begin{equation*}
		\begin{split}
		\langle dE_C(w,\theta),\varphi\rangle=&-H\int_{0}^{L}\bigg(\frac{(w+\ell\sin\theta+y)_{x}\xi(x)}{\sqrt{1+[(w+\ell\sin\theta+y)_{x}]^2}}\bigg)_{x}\varphi\hspace{1mm}dx+\\&-\frac{AE_c}{L_{c}}\Gamma(w+\ell\sin\theta)\int_{0}^{L}\bigg(\frac{(w+\ell\sin\theta+y)_{x}}{\sqrt{1+[(w+\ell\sin\theta+y)_{x}]^2}}\bigg)_{x}\varphi\hspace{1mm}dx;
		\end{split}
		\label{variation w energy Hc}
		\end{equation*}
		by computing the variation of the energy (\ref{energy cavi}) with respect to $\theta$ and integrating by parts, we find for all $\psi\in C^{\infty}_c(0,L)$
		\begin{equation*}
		\begin{split}
		\langle dE_C(w,\theta),\psi\rangle=&-H\ell\int_{0}^{L}\cos\theta\bigg(\frac{(w+\ell\sin\theta+y)_{x}\xi(x)}{\sqrt{1+[(w+\ell\sin\theta+y)_{x}]^2}}\bigg)_{x}\psi\hspace{1mm}dx+\\&-\frac{AE_c\ell}{L_{c}}\Gamma(w+\ell\sin\theta)\int_{0}^{L}\cos\theta\bigg(\frac{(w+\ell\sin\theta+y)_{x}}{\sqrt{1+[(w+\ell\sin\theta+y)_{x}]^2}}\bigg)_{x}\psi\hspace{1mm}dx;
		\end{split}
		\label{variation teta energy Hc}
		\end{equation*}
		similar computations can be performed for the second cable.
	\end{itemize}
	Adding all the energetic contributes of the system, we find   
	\begin{equation}
	\begin{split}
	&\mathcal{E}(t):=\int_{0}^{L}\bigg(\frac{M}{2}w_{t}^2\hspace{1mm}+\frac{M\ell^2}{6}\theta_{t}^2\bigg)\hspace{1mm}dx+\int_{0}^{L}\bigg(\frac{EI}{2}w_{xx}^2+\frac{EJ}{2}\theta_{xx}^2+\frac{GK}{2}\theta_{x}^2\bigg)dx+\\&\hspace{11mm}+H\int_{0}^{L}(\xi\sqrt{1+[(w+\ell\sin\theta+y)_{x}]^2}+\xi\sqrt{1+[(w-\ell\sin\theta+y)_{x}]^2})dx-2H\int_{0}^{L}\xi^2dx+\\&\hspace{11mm}+\frac{AE_c}{2L_{c}}\big([\Gamma(w+\ell\sin\theta)]^2+[\Gamma(w-\ell\sin\theta)]^2\big)-Mg\int_{0}^{L}w\hspace{1mm}dx,
	\label{energy tot}
	\end{split}
	\end{equation}
	that is conserved in time.
	\subsection{The system of evolution partial differential equations}
	\label{The Euler-Lagrange Equations}
	From the energy balance we derive the following system of equations. The unknowns are $w(x,t)$ and $\theta(x,t)$ for $(x,t)\in (0,L)\times(0,\infty)$
	\small
	\begin{equation}
	\begin{cases}
	Mw_{tt} =-EIw_{xxxx}+H\bigg(\dfrac{(w+\ell\sin\theta+y)_{x}\xi}{\sqrt{1+[(w+\ell\sin\theta+y)_{x}]^2}}+\dfrac{(w-\ell\sin\theta+y)_{x}\xi}{\sqrt{1+[(w-\ell\sin\theta+y)_{x}]^2}}\bigg)_{x}+\\\hspace{13mm}+\dfrac{AE_c}{L_{c}}\Gamma(w+\ell\sin\theta)\bigg(\dfrac{(w+\ell\sin\theta+y)_x}{\sqrt{1+[(w+\ell\sin\theta+y)_{x}]^2}}\bigg)_x+\\\vspace{6mm}\hspace{13mm}+\dfrac{AE_c}{L_{c}}\Gamma(w-\ell\sin\theta)\bigg(\dfrac{(w-\ell\sin\theta+y)_{x}}{\sqrt{1+[(w-\ell\sin\theta+y)_{x}]^2}}\bigg)_x+Mg\\
	\frac{M\ell^2}{3}\theta_{tt} =-EJ\theta_{xxxx}+GK\theta_{xx}+H\ell\cos\theta\bigg(\dfrac{(w+\ell\sin\theta+y)_{x}\xi}{\sqrt{1+[(w+\ell\sin\theta+y)_{x}]^2}}-\dfrac{(w-\ell\sin\theta+y)_{x}\xi}{\sqrt{1+[(w-\ell\sin\theta+y)_{x}]^2}}\bigg)_x+\\\hspace{13mm}+\dfrac{AE_c\ell}{L_{c}}\cos\theta\hspace{2mm}\Gamma(w+\ell\sin\theta)\bigg(\dfrac{(w+\ell\sin\theta+y)_{x}}{\sqrt{(1+[w+\ell\sin\theta+y)_{x}]^2}}\bigg)_x+\\\hspace{13mm}-\dfrac{AE_c\ell}{L_{c}}\cos\theta\hspace{2mm}\Gamma(w-\ell\sin\theta)\bigg(\dfrac{(w-\ell\sin\theta+y)_{x}}{\sqrt{1+[(w-\ell\sin\theta+y)_{x}]^2}}\bigg)_x
	\end{cases}
	\label{eq sist}
	\end{equation}
	\normalsize
	where $y(x)$ and $\xi(x)$ depend only on $x$, as defined respectively in (\ref{y(x)})-(\ref{xi}), and $\Gamma(\cdot)$ is defined in (\ref{gamma}); the problem is completed by the boundary and initial conditions:
	\begin{equation}
	\begin{split}
	&w(0,t)=w(L,t)=w_{xx}(0,t)=w_{xx}(L,t)=0 \hspace{7mm}{\rm for } \;t\in (0,\infty)\\&\theta(0,t)=\theta(L,t)=\theta_{xx}(0,t)=\theta_{xx}(L,t)=0 \hspace{11mm}{\rm for } \;t\in (0,\infty)
	\end{split}
	\label{bcs}
	\end{equation}
	\begin{equation}
	\begin{split}
	&w(x,0)=w^0(x),\hspace{4mm} \theta(x,0)=\theta^0(x) \hspace{25mm} {\rm for } \;x\in (0,L)\\
	&w_{t}(x,0)=w^1(x),\hspace{4mm} \theta_{t}(x,0)=\theta^1(x) \hspace{23mm} {\rm for } \;x\in (0,L).\\
	\end{split}
	\label{ics}
	\end{equation}
	We want now to find a weak formulation of (\ref{eq sist}); to do this we consider the Hilbert spaces $L^2(0,L)$, $H^1_{0}(0,L)$, $H^2\cap H^1_{0}(0,L)$ endowed respectively with the scalar products
	\begin{equation*}
	(u,v)_{2}=\int_{0}^{L}uv, \hspace{5mm} (u,v)_{H^1}=\int_{0}^{L}u'v',\hspace{5mm} (u,v)_{H^2}=\int_{0}^{L}u''v''
	\end{equation*}
	and we denote by $H^{-1}(0,L)$ and $H^*(0,L)$ the dual spaces respectively of $H^1_0(0,L)$ and $H^2\cap H^1_{0}(0,L)$ with the corresponding duality $\langle \cdot,\cdot\rangle_1$ and $\langle\cdot,\cdot\rangle_{*}$.\\
	For simplicity we introduce the map $\boldsymbol{\chi}: C^1[0,L]\rightarrow C^0[0,L]$ defined by
	\begin{equation}
	u\mapsto\boldsymbol{\chi}(u):=\dfrac{(u+y)_{x}}{\sqrt{1+[(u+y)_{x}]^2}}.
	\label{chi}
	\end{equation}
	Computing the derivative of $\boldsymbol{\chi}$ with respect to $x$, we obtain the cables curvature along the main span. In the initial configuration ($w=0$), after hooking the deck, the curvature is
	\begin{equation*}
	[\boldsymbol{\chi}(0)]_x=\dfrac{-8f}{L^2\sqrt{\big(1+\frac{64f^2}{L^4}\big(\frac{L}{2}-x\big)^2\big)^3}}\hspace{8mm}\forall x\in(0,L).
	\end{equation*}
	To simplify further the notation we put
	\begin{equation}
	\begin{split}
	&h_\alpha(w,\theta):=-\bigg(H\xi+\dfrac{AE_c}{L_c}\Gamma(w+\ell\sin\theta)\bigg)\boldsymbol{\chi}(w+\ell\sin\theta),\\&h_\beta(w,\theta):=-\bigg(H\xi+\dfrac{AE_c}{L_c}\Gamma(w-\ell\sin\theta)\bigg)\boldsymbol{\chi}(w-\ell\sin\theta),
	\end{split}
	\label{g1,g2}
	\end{equation}
	then (\ref{eq sist}) becomes
	\begin{equation}
	\begin{cases}
	Mw_{tt} =-EIw_{xxxx}-\big[h_\alpha(w,\theta)+h_\beta(w,\theta)\big]_x+Mg\\
	\dfrac{M\ell^2}{3}\theta_{tt} =-EJ\theta_{xxxx}+GK\theta_{xx}-\ell\cos\theta\big[h_\alpha(w,\theta)-h_\beta(w,\theta)\big]_x.
	\end{cases}
	\label{eq sist3}
	\end{equation} 
	with the boundary conditions (\ref{bcs}) and the initial data (\ref{ics}), that we recall here
	\begin{equation}
	\begin{split}
	&w(x,0)=w^0(x),\hspace{4mm} \theta(x,0)=\theta^0(x) \hspace{4mm}\forall x\in (0,L)\\
	&w_{t}(x,0)=w^1(x),\hspace{4mm} \theta_{t}(x,0)=\theta^1(x) \hspace{4mm} \forall x\in (0,L)
	\end{split}
	\label{ics2}
	\end{equation}
	with the regularity 
	\begin{equation}
	w^0, \theta^0\in H^2\cap H^1_0(0,L), \hspace{5mm}w^1,\theta^1\in L^2(0,L).\\
	\label{regularity ics}
	\end{equation}
	We say that $(w,\theta)$ is a weak solution of (\ref{eq sist3}) if $(w, \theta)\in X^2_T$, where
	\begin{equation}
	\begin{split}
	X_T:=C^0\big([0,T];H^2\cap H^1_0(0,L)\big)\cap C^1\big([0,T];L^2(0,L)\big)\cap C^2\big([0,T];H^*(0,L)\big)\\
	\end{split}
	\label{regularity weak}
	\end{equation}
	and if $(w,\theta)$ satisfies the following equations
	\begin{equation}
	\begin{cases}
	M\langle w_{tt},\varphi \rangle_{*}+EI(w,\varphi)_{H^2}=\big(h_\alpha(w,\theta)+h_\beta(w,\theta),\varphi_x\big)_2+\big( Mg,\varphi\big)_2\\
	\dfrac{M\ell^2}{3}\langle\theta_{tt},\psi\rangle_*+EJ(\theta,\psi)_{H^2}+GK(\theta,\psi)_{H^1}=\ell\big(h_\alpha(w,\theta)-h_\beta(w,\theta),(\psi\cos\theta)_x\big)_2
	\end{cases}
	\label{eq weak sist2}
	\end{equation}
	for all $\varphi, \psi\in H^2\cap H^1_0(0,L)$ and $t>0$. Note that in the space $X_T$ the boundary conditions (\ref{bcs}) are already included, then from now we will not mention them. \\In this framework we are ready to state the result about existence and uniqueness of a weak solution
	\begin{teo}
		Let $T>0$ (including the case $T=\infty$), then for all $w^0, \theta^0, w^1,\theta^1$ satisfying (\ref{regularity ics}) there exists a unique (global in time) weak solution $(w,\theta)\in X^2_{T}$ of (\ref{eq sist3}) which satisfies (\ref{ics2}).
		\label{teo1}  
	\end{teo}
	\noindent This result is achieved applying the Galerkin procedure for the existence part and testing the equations with the Green function, applied to the time derivative of the solutions, for the uniqueness part. The presence of the nonlinearities makes challenging the proof that is fully given in Section \ref{Proof of the Theorem}.
	\section{Numerical experiments}
	\label{Numerical experiments}
\begin{figure}[h]
		\centering
		\includegraphics[width=16cm]{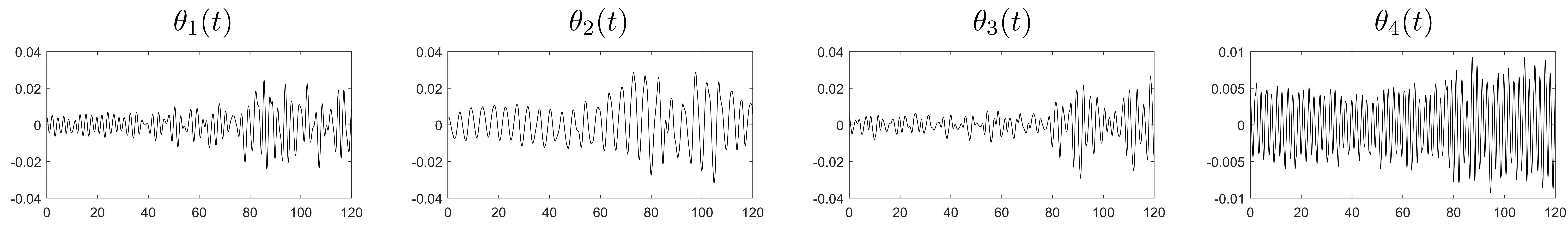}
		\caption{Plots of $\theta_k(t)$ ($k=1,\dots,4$) in radians on $[0,120s]$ with $\overline{w}_9^0=3.87$m.}
		\label{analisi}
	\end{figure}

	In this section we present some numerical experiments on the system (\ref{eq sist})-(\ref{bcs})-(\ref{ics}); in the spirit of the proof of Theorem \ref{teo1} (see Section \ref{Proof of the Theorem}), we apply the Galerkin procedure. More precisely, given the boundary conditions, we seek approximated solutions in the form
	\begin{equation}
	w(x,t)=\sum_{k=1}^{10} w_k(t) \hspace{1mm}e_k, \hspace{5mm} \theta(x,t)=\sum_{k=1}^{4} \theta_k(t) \hspace{1mm}e_k
	\label{approx}
	\end{equation}
	where $e_{k}(x)=\sqrt{\frac{2}{L}}\sin\bigg(\frac{k\pi x}{L}\bigg)$ and  $\sqrt{\frac{2}{L}}$ is a pure number with no unit of measure.
	\begin{defn}
		We call $\overline{w}_k(t):=\sqrt{\frac{2}{L}}w_k(t)$ \textbf{k-th longitudinal mode} and $\overline{\theta}_k(t):=\sqrt{\frac{2}{L}}\theta_k(t)$ \textbf{k-th torsional mode}.
		\label{modi}
	\end{defn}
	\noindent We consider 14 modes because it is a good compromise between limiting computational burden and the possibility of highlighting the instability phenomena which we are interested in. Moreover, from \cite{TNB} we know that, before and during the collapse, the TNB displayed the first 10 longitudinal modes and the second torsional one.\\
	Plugging (\ref{approx}) into (\ref{eq sist}) and projecting onto the space spanned respectively by the first 10 longitudinal modes and the first 4 torsional modes, we obtain an ODE system of 14 equations as (\ref{eq weak sist lin gal  2n}) with the initial conditions 
	\begin{equation*}
	\begin{split}
	&w_k(0)=w^0_k=(w^0,e_k)_2, \hspace{10mm}\dot{w}_k(0)=w^1_k=(w^1,e_k)_2,\hspace{5mm}\forall k=1,\dots,10\\&\theta_k(0)=\theta^0_k=(\theta^0,e_k)_2,\hspace{13mm}
	\dot{\theta}_k(0)=\theta^1_k=(\theta^1,e_k)_2\hspace{9.5mm}\forall k=1,\dots,4.
	\end{split}
	\end{equation*}
	Following Definition \ref{modi}, we put $\overline{w}_k^0:=\sqrt{\frac{2}{L}}w_k^0$, $\overline{w}_k^1:=\sqrt{\frac{2}{L}}w_k^1$ and similarly for the $\theta$ initial conditions.\\
	Applying a similar procedure to \cite{Gazzola Torsion,falocchi} we excite one single longitudinal mode (the $j^{th}$) at a time, applying an initial condition 10$^{-3}$ smaller on all the other components, i.e. in dimensionless form
	\begin{equation}
	\begin{split}
	&\overline{w}_{k}^0=10^{-3}\cdot \overline{w}_{j}^0,\hspace{20mm}\forall k\neq j,\\&\overline{\theta}_{k}^0=\overline{w}_{k}^1=\overline{\theta}_{k}^1=10^{-3}\cdot \overline{w}_{j}^0,\hspace{3mm}\forall k.
	\end{split}
	\label{ic22}
	\end{equation}
	Our aim is to verify if there is a torsional mode that after some time captures energy from the longitudinal modes and becomes larger and larger.\\
	The numerical results are obtained with the MATLAB$^{\textregistered}$ ODE solver ode23tb on the integration time $[0,120s]$, adopting the mechanical constants of the TNB as in Table \ref{tab mec}, see Section \ref{The influence of the mechanical parameters on the stability of the system}; we refer to Section \ref{The influence of the mechanical parameters on the stability of the system} for an analysis of sensitivity in terms of stability of the system by the mechanical parameters.\\
	For each longitudinal mode excited it is possible to determine an instability threshold; but, what do we mean for torsional instability threshold? To give a precise definition in quantitative terms is a hard work. In \cite{garrione} the authors give a definition of instability with respect to the time lapse  considered $[0, T]$ with $T>0$, the amplitude of oscillation and its speed of growth. We think that, choosing a proper interval $[0,T]$ and an appropriate rate of growth of the amplitude, such definition can be useful as a quantitative indicator of instability.\\   
	In our simulations we choose $T=120$s and, following \cite{garrione}, we consider the k-th longitudinal mode unstable if at least one torsional mode grows about 1 order in amplitude in this time lapse; we know that the wide oscillations at the TNB lasted several hours, but we are focusing on the mechanism related to the transfer of energy between longitudinal and torsional modes and this change happened suddenly. Moreover 2 minutes seem to be a sufficient amount of time to consider the system isolated, in which the injection of energy deriving from the wind and the structural capacity to dissipate it, are almost balanced.\\ 
	\begin{figure}[h]
		\centering
		\includegraphics[width=16cm]{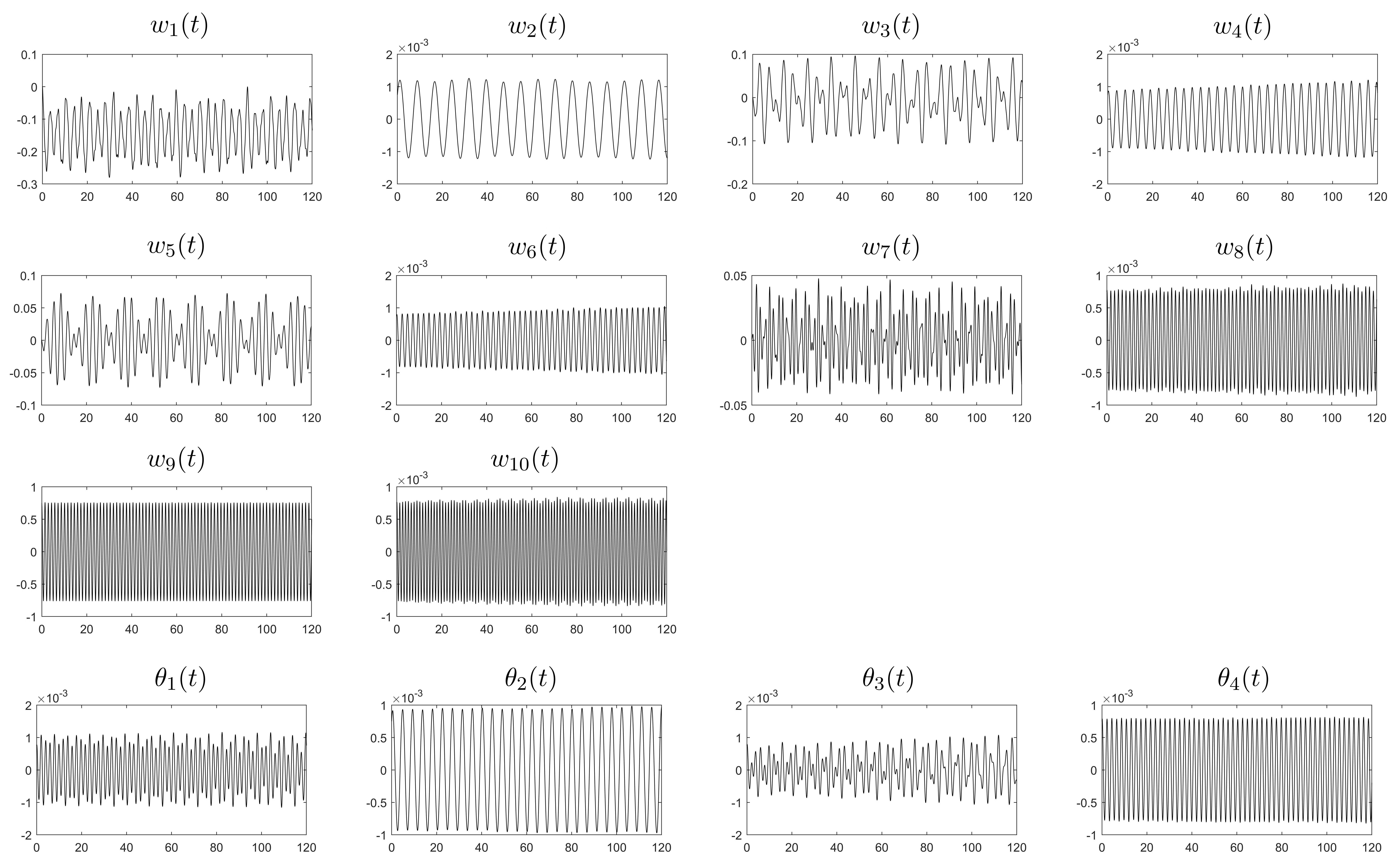}
		\caption{Plots of $w_k(t)$ ($k=1,\dots,10$) in meters and $\theta_k(t)$ ($k=1,\dots,4$) in radians on $[0,120s]$ with $\overline{w}_9^0=0.75$m.}
		\label{analisi stab}
	\end{figure} 
	About the meaning of the instability threshold, we know from the Mathieu equation \cite{mathieu} that the instability with respect to the initial conditions and other parameters can arise in regions called \textit{resonance tongues}, that becomes larger and larger as the initial energy increases; hence, to define a threshold of instability $\overline{w}_{j,th}^0$, as a watershed between stable and unstable situations is not obvious at all. Although our numerical simulations exhibit that for all $j$ there exists $\overline{w}^0_{j,th}>0$ such that for $\overline{w}_j^0<\overline{w}_{j,th}^0$ the torsional components are stable, while for $\overline{w}_j^0>\overline{w}_{j,th}^0$ they are unstable, it may happen that there are very thin \textit{resonance tongues} for some initial conditions below the threshold of instability $\overline{w}_{j,th}^0$. We neglect these cases since the probability to fall in a thin \textit{resonance tongue} is small, and, even if we were in this situation, it is very probable that the solutions come back to a stability region in a while.\\
	In this section we focus our attention on high longitudinal modes and low torsional modes; indeed, from \cite[p.29]{TNB} we know that the morning of the TNB failure 
	\begin{quote}
		the center span was oscillating with either 8 or 9 nodes [i.e. as $\sin(\frac{9\pi x}{L})$ and  $\sin(\frac{10\pi x}{L})$].  [$\dots$], at 10:00 A.M. the center span developed a torsional movement with a node at mid span [i.e. as $\sin(\frac{2\pi x}{L})$].
	\end{quote}
	In Figures \ref{analisi stab} and \ref{analisi} we report the results of two analysis on the system, imposing respectively $\overline{w}_9^0=$0.75m and $\overline{w}_{9}^0=$3.87m. As we can see, Figure \ref{analisi stab} presents a situation of stability, while in Figure \ref{analisi}, where we show only the torsional modes for brevity, we are close to the torsional instability threshold and the first 3 torsional modes after 80s suddenly begin to grow.
	These results reveal that there is an exchange of energy between longitudinal and torsional modes, due exclusively to the initial amplitude (and then, to the initial energy) of the longitudinal mode excited. Moreover we point out that the oscillations of the 2$^{nd}$ torsional mode increase quantitatively more than the others, confirming the real observations on TNB collapse.
	Our results show that the 9$^{th}$ and 10$^{th}$ longitudinal modes are very prone to develop torsional instability; in these cases we record $\overline{w}_{9,th}^0\approx 3.87$m and $\overline{w}_{10,th}^0\approx3.40$m lower with respect, for instance, the thresholds of 7$^{th}$ and 8$^{th}$ modes, where $\overline{w}_{7,th}^0\approx 4.90$m and $\overline{w}_{8,th}^0\approx5.15$m. 
	\section{The influence of the mechanical parameters on the stability of the system}
	\label{The influence of the mechanical parameters on the stability of the system}
\begin{figure}[h]
		\centering
		\includegraphics[width=16cm]{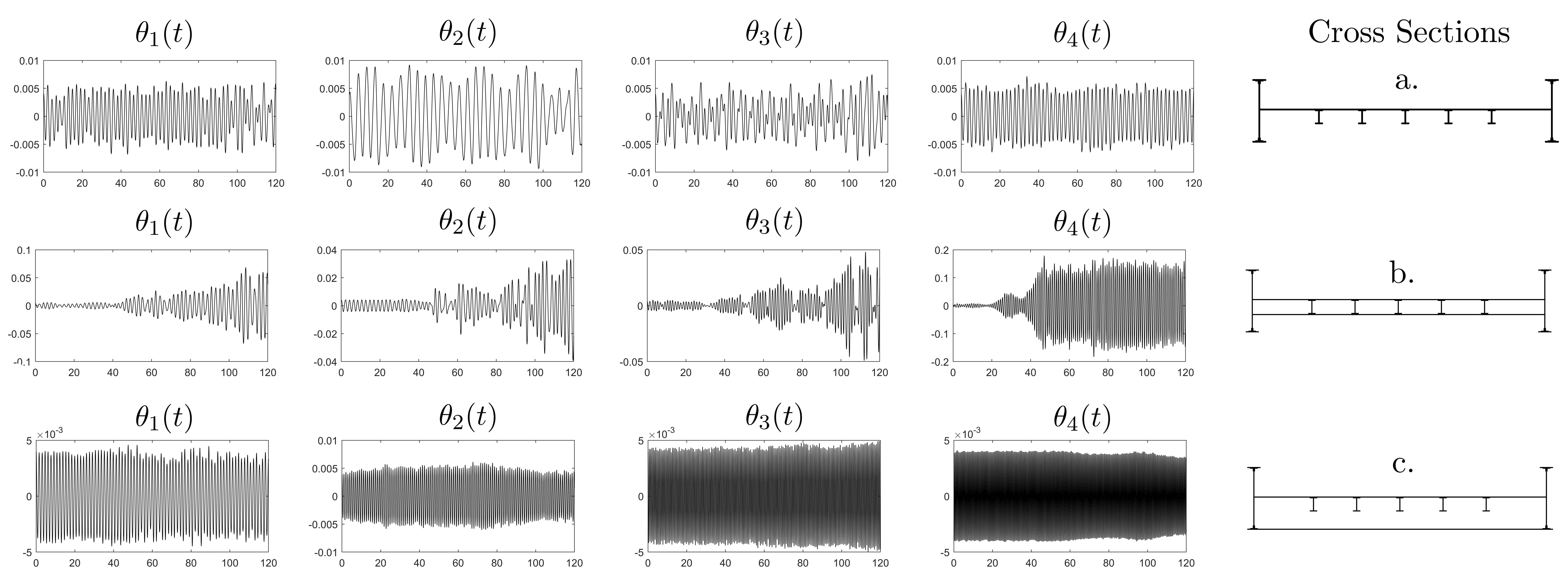}
		\caption{Comparison between the torsional modes in the case a. $I=1.54$m$^4$, b. $K=0.1337$m$^4$ and c. $K=0.7171$m$^4$ with the corresponding cross sections.}
		\label{EI}
	\end{figure}

	The system (\ref{eq sist}) depends on several mechanical constants that characterize the suspension bridge. In this section we study how the torsional instability of this system is affected by these parameters.\\ As in Section \ref{Numerical experiments}, we excite the 9$^{th}$ longitudinal mode ($\overline{w}_9^0=3.87$m), applying an initial condition 10$^{-3}$ smaller on all the others components on $[0,120s]$.
	We are interested more on the qualitative datum respect to the quantitative; for brevity we do not show the plots of the 10 longitudinal modes.\\
	We denote as "basic situation" the solution of the system with the mechanical properties of the TNB, listed in Table \ref{tab mec} (values taken from \cite{TNB,Plaut}). Note that the constants $H$ and $L_c$ depend on the previous by the equations (\ref{H}) and (\ref{gamma}). Hence, in our model the behavior of the suspension bridge is influenced by 11 parameters; among them there are standard values in the bridge design while others highly depend on the designer choice. Typically, when a bridge is built, the length of the main span is fixed with respect to the site conditions, and, consequently, the width of the roadway; for these reasons in our numerical experiments we maintain fixed the values $L$ and $\ell$.\\
	The usual material employed to build the bearing structure is the steel and then we consider quite reliable $E$ and $G$, the Young and shear modulus of the deck; on the other hand, the elastic modulus of the cables has to be reduced with respect to percentage of air void and the kind of ropes used in the assemblage. In \cite{Plaut} $E_c=185$GPa is considered the conventional value in the design of suspension bridges, moreover, other values of $E_c$, defined in \cite{Podolny} for every kind of ropes, remain quite close to the previous. For these reasons we do not modify the elastic constants.\\
	\begin{table}[h]
		\begin{center}
			\begin{tabular}{*{3}{c}}
				\hline
				$E$:&210\hspace{1mm}000MPa& Young modulus of the deck (steel)\\
				$E_c$:&185\hspace{1mm}000MPa& Young modulus of the cables (steel)\\
				$G$:&81\hspace{1mm}000MPa& Shear modulus of the deck (steel)\\
				$L$:&853.44m& Length of the main span\\
				$\ell:$&6m& Half width of the deck\\
				$f$:&70.71m& Sag of the cable\\
				$I$:&0.154m$^4$& Moment of inertia of the deck cross section\\
				$K$:&6.07$\cdot10^{-6}$m$^4$& Torsional constant of the deck\\
				$J$:&5.44m$^6$& Warping constant of the deck\\
				$A$:&0.1228m$^2$& Area of the cables section\\
				$M$:&7198kg/m&Mass linear density of the deck \\
				$H$:&45\hspace{1mm}413kN& Initial tension in the cables, see (\ref{H})\\
				$L_c$:&868.815m& Initial length of the cables, see (\ref{gamma})\\
				\hline
			\end{tabular}
		\end{center}
		\caption{TNB mechanical features.}
		\label{tab mec}
	\end{table}
	The sag-span ratio $\frac{f}{L}$ assumes an important role in the bridge behavior, affecting the horizontal component of the cable force $H$ and the total stiffness of the bridge; in the design practice
	\begin{equation*}
	\dfrac{f}{L}=\dfrac{1}{12}\div\dfrac{1}{8}
	\end{equation*}
	and more the ratio is large more the stresses are minimized \cite{Podolny}. In the TNB $\frac{f}{L}\approx \frac{1}{12}$, probably due to the requirement to reduce the tower height in order to have an economic safe.\\
	\begin{figure}[h]
		\centering
		\includegraphics[width=16cm]{"instabilityrigid3"}
		\includegraphics[width=16cm]{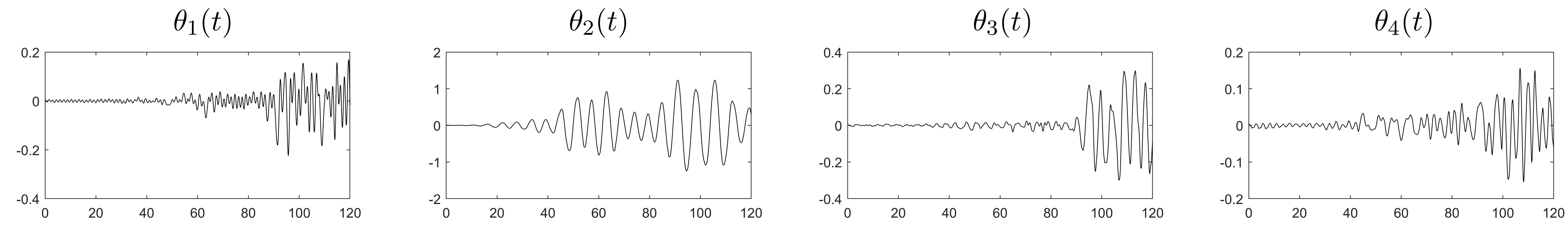}
		\caption{Comparison between the torsional modes in the case $f=70.71$m, i.e. the basic situation, (above) and $f=106.71$m (below).}
		\label{f}
	\end{figure}
	From (\ref{eq sist}) we observe that the sag-span ratio is highly involved in the system, e.g. in $H$, $L_c$, $y'(x)$, $\xi(x)$ and $\Gamma(\cdot)$. In Figure \ref{f} we compare the first 4 torsional modes in the cases $f=70.71$m and $f=106.71$m, which correspond respectively to a sag-span ratio equal to $\frac{1}{12}$ and $\frac{1}{8}$. It turns out that 
	\begin{quote}
		\textbf{an increment of the sag $f$ determines a larger torsional instability in the bridge.}
	\end{quote}
	In particular, it is interesting to note that when $f=106.71$m there is a 30\% decrement of $H$ towards a 2\% decrement of $\frac{AE_c}{L_c}$. Then we have that the torsional instability of the system is sensitive to the constants $H$ and $\frac{AE_c}{L_c}$ and it grows when $H$ decreases and $\frac{AE_c}{L_c}$ increases. \\A further confirmation of it appears if we increase the sectional area of the cable $A$; note that for static reasons is not possible to reduce too much $A$ and on the other hand for practical reasons (installation and tensioning) to increase overly. In any case a designer should look for reducing the sectional area of the cable not only for the stability aim deriving from our model, but also because, as the cable becomes so heavy, its capability to carry live load decreases, as suggested by \cite[p.30]{Podolny}.\\
	The torsional stability of the system can be improved also modifying the geometry of the deck's section. In general, the torsional performance of closed cross sections
	is better than that of open sections \cite{usainstitute}; the cross section of the TNB was open and this is one of the reasons why it was very prone to develop torsional instability. Indeed, after its failure most long span bridges were built with closed cross section increasing their stiffness (truss-stiffened section). In Figure \ref{EI} are plotted the $\theta_k(t)$ ($k=1,\dots,4$) components varying the moment of inertia $I$ (case a.) and the torsional constant of the deck $K$ (case b. and c.).
	With respect to the basic situation we observe that increasing of 1 order $I$, e.g. enhancing the thickness of the profiles maintaining an open cross section, we gain in terms of torsional stability. Also acting on the torsional properties of the section can give good results. For instance, comparing Figure \ref{EI}b. and \ref{EI}c. we note that the introduction of a continuous plate of 2cm in the correct position (on the bottom) reduces considerably the torsional instability and produces a relevant growth in the $K$ constant; in fact, only closing the cross section, $K$ rises of 5 orders of magnitude! Then our model shows that
	\begin{quote}
		\textbf{a deck with closed cross section is torsionally more stable than the same deck with open section. }
	\end{quote}
	About the warping constant $J$ we record that, modifying the section properties in a physical way, it does not change enough to be considered significant in terms of torsional stability of the bridge.\\
	Last but not least, the linear mass $M$ of the deck is another important parameter to prevent the torsional instability, indeed, an increased mass determines a greater energy-storage capacity of the structure, reducing the oscillation's amplitudes \cite{Podolny}. In (\ref{eq sist}) $M$ is involved in the inertia terms and implicitly in the constant $H$; to enhance $M$ implies an increment of $H$ and, as discussed before, this fact acts in favor of stability. From a database on suspension bridges published in \cite[p.91]{Podolny}, normalizing the masses to the bridges width, it turns out that TNB had a linear mass approximately equal to $40\%\div60\%$ the linear mass of the other bridges; even if the others have a span $20\%\div30\%$ longer, the datum on the TNB mass is surprising and gives a further justification on the torsional oscillations recorded during its collapse. 
	\begin{figure}[h]
		\centering
		\includegraphics[width=16cm]{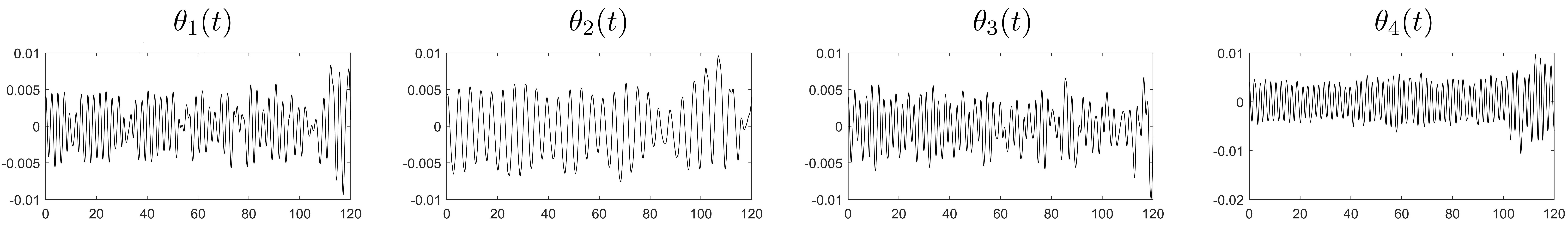}
		\caption{The torsional modes for the TNB with $M=10077\frac{kg}{m}$.}
		\label{M}
	\end{figure}  
	In Figure \ref{M} we show $\theta_k(t)$ ($k=1,\dots,4$) with the original linear mass of TNB 40\% increased; although we are under the average normalized linear mass of the other bridges in \cite[p.91]{Podolny}, the results is meaningful, because we see how
	\begin{quote}
		\textbf{an increment in the mass of the deck reduces the torsional instability.}
	\end{quote}
	In this section we have considered separately all the parameters that play a role in suspension bridges stability; we point out that the optimal situation in structural and also in economical sense can be achieved with an accurate combination of all these parameters. For instance, the increment of $M$ usually implies an increment of $I$ and $K$, since the cross section is heavier (and thicker), in this way the best solution is not necessarily the most expensive.  
	
	\section{Proof of the Theorem \ref{teo1}}
	\label{Proof of the Theorem}
	In this Section we prove the existence and uniqueness of a weak solution of (\ref{eq sist3}). The proof uses a Galerkin procedure and is divided in several steps as in classical hyperbolic PDE problems, see for instance \cite{ball,Berkovits,Gazzola hyperb, Holubova}.
	\subsection{Existence of solutions}
	\label{Existence of solutions}
	\textbf{Step 1: Construction of a sequence of solutions in finite dimensional spaces approximating $\boldsymbol{(w,\theta)}$} \\
	An orthogonal basis of $L^2(0,L)$, $H^1_{0}(0,L)$, $H^2\cap H^1_{0}(0,L)$  is $\{e_{k}\}^{\infty}_{k=1}$, where 
	\begin{equation*}
	e_{k}(x)=\sqrt{\frac{2}{L}}\sin\bigg(\frac{k\pi x}{L}\bigg), \hspace{5mm}|| e_{k}||_{2}=1,\hspace{5mm}||e_{k}||_{H^1}=\frac{k\pi}{L},\hspace{5mm}|| e_{k}||_{H^2}=\frac{k^2\pi^2}{L^2},
	\end{equation*}
	then for any $n\geq 1$ we introduce the space
	\begin{equation*}
	E_n:={\rm span} \{e_1,\dots,e_n\}.
	\end{equation*}
	We put for any $n\geq 1$
	\begin{equation*}
	\begin{split}
	&w^0_n:=\sum_{k=1}^{n} (w^0,e_k)_2 \hspace{1mm}e_k =\frac{L^4}{\pi^4}\sum_{k=1}^{n}\dfrac{(w^0,e_k)_{H^2}}{k^4} \hspace{1mm}e_k, \\
	&\theta^0_n:=\sum_{k=1}^{n} (\theta^0,e_k)_2 \hspace{1mm}e_k =\sum_{k=1}^{n}\bigg(EJ\frac{k^4\pi^4}{L^4}+GK\frac{k^2\pi^2}{L^2}\bigg)^{-1}[EJ(\theta^0,e_k)_{H^2}+GK(\theta^0,e_k)_{H^1}]\hspace{1mm}e_k,\\& w^1_n:=\sum_{k=1}^{n} (w^1,e_k)_2 \hspace{1mm}e_k,  \hspace{8mm} \theta^1_n:=\sum_{k=1}^{n}(\theta^1,e_k)_2 \hspace{1mm}e_k, 
	\end{split}
	\end{equation*}
	so that
	\begin{equation}
	\begin{split}
	&w^0_n\rightarrow w^0 \hspace{2mm}{\rm in}\hspace{2mm} H^2, \hspace{5mm} \theta^0_n\rightarrow\theta^0 \hspace{2mm}{\rm in}\hspace{2mm} H^2,\hspace{5mm}
	w^1_n\rightarrow w^1 \hspace{2mm}{\rm in}\hspace{2mm} L^2, \hspace{5mm}\theta^1_n\rightarrow\theta^1 \hspace{2mm}{\rm in}\hspace{2mm} L^2
	\end{split}
	\label{converg ics}
	\end{equation}
	as $n\rightarrow\infty$. For any $n\geq 1$ we seek $(w_n,\theta_n)$ such that
	\begin{equation*}
	w_n(x,t)=\sum_{k=1}^{n} w_n^k(t) \hspace{1mm}e_k, \hspace{5mm} \theta_n(x,t)=\sum_{k=1}^{n} \theta_n^k(t) \hspace{1mm}e_k
	\end{equation*}
	and which solves the problem (\ref{eq weak sist2}).
	Using as test functions $\varphi, \psi\in E_n$, (\ref{eq weak sist2}) becomes
	\begin{equation}
	\begin{cases}
	M\big((w_n)_{tt},e_j\big)_2+EI\big(w_n,e_j\big)_{H^2}=\big(h_\alpha(w_n,\theta_n)+h_\beta(w_n,\theta_n),e'_j\big)_2+(Mg,e_j)_2\\
	\frac{M\ell^2}{3}\big((\theta_n)_{tt},e_j\big)_2+EJ\big(\theta_n,e_j\big)_{H^2}+GK\big(\theta_n,e_j\big)_{H^1}=\ell\big(h_\alpha(w_n,\theta_n)-h_\beta(w_n,\theta_n),(e_j\cos\theta_{n})_x\big)_2.
	\end{cases}
	\label{eq weak sist gal}
	\end{equation}
	Testing $n$ times each equation for $j=1,\dots,n$ we obtain 
	\begin{equation}
	\begin{cases}
	M\ddot{w}^k_n(t)+EI\dfrac{k^4\pi^4}{L^4} w^k_n(t)=\big(h_\alpha(w_n,\theta_n)+h_\beta(w_n,\theta_n),e'_k\big)_2+Mg\dfrac{\sqrt{2L}(1-(-1)^k)}{k\pi}\\
	\dfrac{M\ell^2}{3}\ddot{\theta}^k_n(t)+\bigg(EJ\dfrac{k^4\pi^4}{L^4}+GK\dfrac{k^2\pi^2}{L^2}\bigg)\theta^k_n(t)=\ell\big(h_\alpha(w_n,\theta_n)-h_\beta(w_n,\theta_n),(e_k\cos\theta_{n})_x\big)_2
	\end{cases}
	\label{eq weak sist lin gal  2n}
	\end{equation}
	\begin{flushright}
		$\forall k=1,\dots,n$.
	\end{flushright}
	Since $h_\alpha(w_n,\theta_n)$ and  $h_\beta(w_n,\theta_n)$ are continuous, from the theory of ODEs this finite-dimensional system with the initial conditions
	\begin{equation*}
	\begin{split}
	&w^k_n(0)=(w^0,e_k)_2, \hspace{10mm}\theta^k_n(0)=(\theta^0,e_k)_2\hspace{10mm}
	\dot{w}^k_n(0)=(w^1,e_k)_2,\hspace{10mm} \dot{\theta}^k_n(0)=(\theta^1,e_k)_2
	\end{split}
	\end{equation*}
	admits a local solution defined on some  $[0,t_n)$ with $t_n\in(0,T]$.\\
	
	\textbf{Step 2: Uniform bounds for the sequence $\boldsymbol{\{(w_n,\theta_n)\}}$} \\
	We omit for the moment the spatial dependence of the approximated solutions. We test the first equation in (\ref{eq weak sist gal}) by $\dot{w}_n$, the second by $\dot{\theta}_n$, we integrate in $x$ on $(0,L)$ and we sum the two equations; then we find
	\begin{equation}
	\begin{split}
	&\dfrac{M}{2}\frac{d}{dt}||\dot{w}_n||^2_2+\dfrac{EI}{2}\frac{d}{dt}||w_n||^2_{H^2}+\dfrac{M\ell^2}{6}\frac{d}{dt}||\dot{\theta}_n||^2_2+\dfrac{EJ}{2}\frac{d}{dt}||\theta_n||^2_{H^2}+\dfrac{GK}{2}\frac{d}{dt}||\theta_n||^2_{H^1}=\\&=-H\int_{0}^{L}\xi\dfrac{(w_n+\ell\sin\theta_n+y)_x}{\sqrt{1+[(w_n+\ell\sin\theta_n+y)_x]^2}}(\dot{w}_n+\ell\dot{\theta}_n\cos\theta_n)_{x}dx+\\&\hspace{4mm}-H\int_{0}^{L}\xi\dfrac{(w_n-\ell\sin\theta_n+y)_x}{\sqrt{1+[(w_n-\ell\sin\theta_n+y)_x]^2}}(\dot{w}_n-\ell\dot{\theta}_n\cos\theta_{n})_{x}dx+\\&\hspace{4mm}-\dfrac{AE_c}{L_c}\int_{0}^{L}\Gamma(w_n+\ell\sin\theta_n)\dfrac{(w_n+\ell\sin\theta_n+y)_x}{\sqrt{1+[(w_n+\ell\sin\theta_n+y)_x]^2}}(\dot{w}_n+\ell\dot{\theta}_n\cos\theta_n)_{x}dx+\\&\hspace{4mm}-\dfrac{AE_c}{L_c}\int_{0}^{L}\Gamma(w_n-\ell\sin\theta_n)\dfrac{(w_n-\ell\sin\theta_n+y)_x}{\sqrt{1+[(w_n-\ell\sin\theta_n+y)_x]^2}}(\dot{w}_n-\ell\dot{\theta}_n\cos\theta_{n})_{x}dx+\\&\hspace{5mm}+\int_{0}^{L}Mg\dot{w}_ndx
	\end{split}
	\label{step 2}
	\end{equation}
	Recalling the energy (\ref{energy tot}), we write (\ref{step 2}) as
	\begin{equation}
	\begin{split}
	&\frac{d}{dt}\bigg[\dfrac{M}{2}||\dot{w}_n||^2_2+\dfrac{EI}{2}||w_n||^2_{H^2}+\dfrac{M\ell^2}{6}||\dot{\theta}_n||^2_2+\dfrac{EJ}{2}||\theta_n||^2_{H^2}+\dfrac{GK}{2}||\theta_n||^2_{H^1}-Mg\int_{0}^{L}w_ndx+\\&+H\int_{0}^{L}(\xi\sqrt{1+[(w_n+\ell\sin\theta_n+y)_x]^2}+\xi\sqrt{1+[(w_n-\ell\sin\theta_n+y)_x]^2})dx+\\&+\dfrac{AE_c}{2L_c}\big([\Gamma(w_n+\ell\sin\theta_n)]^2+[\Gamma(w_n-\ell\sin\theta_n)]^2\big)\bigg]=0.
	\end{split}
	\label{step 2.1}
	\end{equation} 
	We denote by 
	\begin{equation}
	\begin{split}
	&\mathcal{E}_n(t):=\int_{0}^{L}\bigg(\frac{M}{2}\dot{w}_n^2+\frac{EI}{2}[(w_n)_{xx}]^2+\frac{M\ell^2}{6}\dot{\theta}_n^2+\frac{EJ}{2}[(\theta_n)_{xx}]^2+\frac{GK}{2}[(\theta_n)_{x}]^2\bigg)dx+\\&\hspace{11mm}+H\int_{0}^{L}(\xi\sqrt{1+[(w_n+\ell\sin\theta_n+y)_{x}]^2}+\xi\sqrt{1+[(w_n-\ell\sin\theta_n+y)_{x}]^2})dx+\\&\hspace{11mm}+\frac{AE_c}{2L_{c}}\big([\Gamma(w_n+\ell\sin\theta_n)]^2+[\Gamma(w_n-\ell\sin\theta_n)]^2\big)-Mg\int_{0}^{L}w_n\hspace{1mm}dx,
	\end{split}
	\label{Fn}
	\end{equation}
	the energy $\mathcal{E}(t)$ of the approximated solution $(w_n,\theta_n)$, introduced in (\ref{energy tot}) up to an additive constant. Integrating (\ref{step 2.1}) in $s$ on $(0,t)$ for $t\in(0,T)$ we find
	\begin{equation}
	\begin{split}
	&\mathcal{E}_n(t)=c_n
	\end{split}
	\label{step 2.2}
	\end{equation} 
	where 
	\small
	\begin{equation*}
	\begin{split}
	&c_n:=\frac{M}{2}||w^1_n||^2_2+\frac{EI}{2}||w^0_n||^2_{H^2}+\frac{M\ell^2}{6}||\theta^1_n||^2_2+\dfrac{EJ}{2}||\theta^0_n||^2_{H^2}+\frac{GK}{2}||\theta^0_n||^2_{H^1}-Mg\int_{0}^{L}w^0_ndx+\\&+H\int_{0}^{L}(\xi\sqrt{1+[(w^0_n+\ell\sin\theta^0_n+y)']^2}+\xi\sqrt{1+[(w^0_n-\ell\sin\theta^0_n+y)']^2})dx+\dfrac{AE_c}{2L_c}\hspace{1mm}\cdot\\&\cdot\bigg\{\bigg(\int_{0}^{L}[\sqrt{1+[(w^0_n+\ell\sin\theta^0_n+y)']^2}-\xi] dx\bigg)^2+\bigg(\int_{0}^{L}[\sqrt{1+[(w^0_n-\ell\sin\theta^0_n+y)']^2}-\xi] dx\bigg)^2\bigg\}.
	\end{split}
	\label{step 2.3}
	\end{equation*} 
	\normalsize
	In (\ref{Fn}) there is only a term with undefined sign, for which there exists $\varepsilon>0$ such that
	\begin{equation*}
	-\int_{0}^{L}w_ndx\geq -C_1\int_{0}^{L}(1+\varepsilon w_n^2)dx\geq-C_1(L+\varepsilon||w_n||_2^2)\geq-C_1(L+C_2\varepsilon||w_n||_{H^2}^2)
	\end{equation*}
	for $C_1, C_2>0$. Since $EI>2MgC_1C_2\varepsilon$, we find $\eta>0$ such that
	\begin{equation*}
	\begin{split}
	&\mathcal{E}_n(t)\geq\frac{M}{2}||\dot{w}_n||_2^2+\bigg(\frac{EI}{2}-MgC_1C_2\varepsilon\bigg)||w_n||_{H^2}^2+\frac{M\ell^2}{6}||\dot{\theta}_n||^2_2+\frac{EJ}{2}||\theta_n||_{H^2}^2+\\&\hspace{10mm}+\frac{GK}{2}||\theta_n||_{H^1}^2+H\int_{0}^{L}(\xi\sqrt{1+[(w_n+\ell\theta_n+y)_{x}]^2}+\xi\sqrt{1+[(w_n-\ell\theta_n+y)_{x}]^2})dx+\\&\hspace{10mm}+\frac{AE}{2L_{c}}\big([\Gamma(w_n+\ell\theta_n)]^2+[\Gamma(w_n-\ell\theta_n)]^2\big)-MgC_1L\geq\\&\hspace{10mm}\geq\eta(||\dot{w}_n||^2_2+||w_n||^2_{H^2}+||\dot{\theta}_n||^2_2+||\theta_n||^2_{H^2}+||\theta_n||^2_{H^1})-MgC_1L.
	\end{split}
	\label{step 2.6}
	\end{equation*}
	Then (\ref{step 2.2}) becomes
	\begin{equation*}
	\begin{split}
	&\eta(||\dot{w}_n||^2_2+||w_n||^2_{H^2}+||\dot{\theta}_n||^2_2+||\theta_n||^2_{H^2}+||\theta_n||^2_{H^1})\leq C_0+MgC_1L.
	\label{step 2.bo}
	\end{split}
	\end{equation*} 
	where the constant $C_0:=\sup_n(|c_n|)<\infty$ is independent on $n$ and finite thanks to (\ref{regularity ics}).
	Then, we have the bound on $(w_n,\theta_n)$
	\begin{equation}
	||\dot{w}_n||^2_2+||w_n||^2_{H^2}+||\dot{\theta}_n||^2_2+||\theta_n||^2_{H^2}+||\theta_n||^2_{H^1}\leq C_3 \hspace{10mm} \forall t\in [0,T].
	\label{step2.8}
	\end{equation}
	Since $C_3$ does not depend on $n$ and $t$, the global existence of $(w_n,\theta_n)$ on $[0,T]$ is ensured.\\
	
	\textbf{Step 3: A strongly convergent subsequence for $\boldsymbol{\{(w_n,\theta_n)\}}$} \\
	To simplify the notation we denote by $L^p(V)$ the space $L^p((0,T);V(0,L))$ for $1\leq p\leq\infty$ and by $Q=(0,T)\times(0,L)$. From the estimate (\ref{step2.8}) we see that
	\begin{equation*}
	\begin{split}
	&\{w_n\}, \{\theta_n\} \hspace{18mm} {\rm are\hspace{2mm} bounded \hspace{2mm} in} \hspace{10mm} L^{\infty}(H^2),\\&\{\dot{w}_n\}, \{\dot{\theta}_n\} \hspace{18mm} {\rm are\hspace{2mm} bounded \hspace{2mm} in} \hspace{10mm} L^{\infty}(L^2).
	\end{split}
	\end{equation*}
	Then, it is possible to extract a subsequence, still denoted by $n$, such that 
	\begin{equation*}
	\begin{split}
	&w_{n}  \stackrel{*}{\rightharpoonup}w, \hspace{3mm} \theta_{n}\stackrel{*}{\rightharpoonup} \theta \hspace{18mm} {\rm in} \hspace{10mm} L^{\infty}( H^2) ,\\&\dot{w}_n\stackrel{*}{\rightharpoonup}z,\hspace{3mm} \dot{\theta}_n\stackrel{*}{\rightharpoonup}\alpha \hspace{18.5mm} {\rm in} \hspace{10mm} L^{\infty}(L^2),
	\end{split}
	\end{equation*}
	in which the symbol $\stackrel{*}{\rightharpoonup}$ indicates the weak* convergence in $L^{\infty}$; from the definition of weak* convergence and distributional derivative we obtain that $\dot{w}=z$ and $\dot{\theta}=\alpha$.\\In particular from the boundedness of $\{w_n\}, \{\theta_n\}$ and $\{\dot{w}_n\}, \{\dot{\theta}_n\}$ we also have weak converge respectively in $L^{2}(H^2)$ and $L^{2}(Q)$; then, due to the compact embedding $H^1(Q)\subset L^2(Q)$, we obtain the strong convergence
	\begin{equation*}
	w_{n}  \rightarrow w, \hspace{3mm} \theta_{n}\rightarrow\theta \hspace{12mm} {\rm in} \hspace{10mm} L^{2}(Q),
	\end{equation*}
	from which $\sin\theta_{n}\rightarrow\sin\theta $ in $L^2(Q)$, since $||\sin\theta_n-\sin\theta||_{L^2(Q)}\leq||\theta_n-\theta||_{L^2(Q)}\rightarrow0$ as $n\rightarrow\infty$, (similarly $\cos\theta_{n}\rightarrow\cos\theta$).\\About the nonlocal term $\Gamma$, defined in (\ref{gamma}), we see that
	\begin{equation*}
	\Gamma(w_n\pm\ell\sin\theta_n)=\int_{0}^{L}\big(\sqrt{1+[(w_n\pm\ell\sin\theta_{n}+y)_x]^2}-\sqrt{1+[y']^2}\big)dx\rightarrow 	\Gamma(w\pm\ell\sin\theta),
	\end{equation*}
	thanks to the Lebesgue's dominated convergence Theorem.\\
	Let now consider the functional $\boldsymbol{\chi}$, defined in (\ref{chi}) and let note that $|\boldsymbol{\chi}(u)|<1$ for all $u\in C^1[0,L]$; then we have that $ \boldsymbol{\chi}^2(w_n\pm\ell\sin\theta_n)< 1$ and
	\begin{equation*}
	||\boldsymbol{\chi}(w_n\pm\ell\sin\theta_n)||^2_{L^2(Q)}=\int_{0}^{T}\int_{0}^{L}\dfrac{[(w_n\pm\ell\sin\theta_n+y)_{x}]^2}{1+[(w_n\pm\ell\sin\theta_n+y)_{x}]^2}dx dt< LT.
	\end{equation*}
	Hence $\boldsymbol{\chi}(w_n\pm\ell\sin\theta_n)$ converges weakly, up to a subsequence, to $\boldsymbol{\chi}(w\pm\ell\sin\theta)$ in $L^2(Q)$ and it is possible to pass to the limit the first equation in (\ref{eq weak sist gal}). \\To do the same for second equation in (\ref{eq weak sist gal}) we consider that 
	\begin{equation*}
	||\boldsymbol{\chi}(w_n\pm\ell\sin\theta_n)\cos\theta_n||^2_{L^2(Q)}< LT\hspace{7mm}||\boldsymbol{\chi}(w_n\pm\ell\sin\theta_n)\theta_{nx}\sin\theta_n||^2_{L^2(Q)}\leq C_7||\theta_n||^2_{L^\infty(H^1)},
	\end{equation*}
	which implies the weak convergence of this terms in $L^2(Q)$.
	Next, recalling the convergence of the initial conditions (\ref{converg ics}), we find that $(w,\theta)$ is a weak solution of (\ref{eq sist3})-(\ref{ics2}), such that $w, \theta\in L^\infty((0,T);H^2\cap H_0^1(0,L))$ and $\dot{w}, \dot{\theta}\in L^\infty((0,T);L^2(0,L))$.\\
	Thanks to Lemma 3.2 \cite[p.69]{Temam} we infer that the components $w, \theta\in C^0([0,T]; L^2(0,L))$ and $\dot{w},\dot{\theta}\in C^0([0,T];H^*(0,L))$. Hence, exploiting this fact and the boundedness of $\{w(t)\}$ and $\{\theta(t)\}$ in $H^2$ (resp. $\{\dot{w}(t)\}$ and $\{\dot{\theta}(t)\}$ in $L^2$), we deduce the weak continuity of the solution respect to time.\\
	The strong continuity can be inferred integrating the energy equality (\ref{step 2.1}) satisfied by $(w,\theta)$, from $(0,t_n)$ and from $(0,t_0)$, subtracting the two results and passing to the limit for all $t_n\rightarrow t_0$. \\Adding from (\ref{eq sist3}) the regularity $w\in C^2([0,T]; H^*(0,L))$, we have proved the existence of a weak solution $(w, \theta)\in X^2_T$ of (\ref{eq sist3}) over the interval $(0,T)$, satisfying (\ref{ics2}); we know that the total energy of (\ref{eq sist}) is conserved in time, then the solution cannot blow up in finite time and the global existence is obtained for an arbitrary $T>0$. 
	\subsection{Uniqueness of the solution}
	\label{Existence and uniqueness result}
	For contradiction, consider two solutions $(w_1,\theta_1)$, $(w_2,\theta_2)\in X^2_T$ satisfying the same initial conditions (\ref{ics2}). By subtracting the two systems satisfied by $(w_j,\theta_j)$ with $j=1,2$ and denoting by
	$w=w_1-w_2$ and $\theta=\theta_1-\theta_2$, we see that $(w,\theta)$ is a solution of
	\begin{equation}
	\begin{cases}
	&M\langle w_{tt},\varphi \rangle_{*}+EI(w,\varphi)_{H^2}=\big(h_\alpha(w_1,\theta_1)-h_\alpha(w_2,\theta_{2}),\varphi_x\big)_2+\big(h_\beta(w_1,\theta_1)-h_\beta(w_2,\theta_{2}),\varphi_x\big)_2\\
	&\dfrac{M\ell^2}{3}\langle\theta_{tt},\psi\rangle_*+EJ(\theta,\psi)_{H^2}+GK(\theta,\psi)_{H^1}=\ell\big(h_\alpha(w_1,\theta_1),(\psi\cos\theta_1)_x\big)_2+\\&-\ell\big(h_\alpha(w_2,\theta_2),(\psi\cos\theta_2)_x\big)_2-\ell\big(h_\beta(w_1,\theta_1),(\psi\cos\theta_1)_x\big)_2+\ell\big(h_\beta(w_2,\theta_2),(\psi\cos\theta_2)_x\big)_2
	\end{cases}
	\label{eq weak sist2sol}
	\end{equation}
	for all $\varphi, \psi\in H^2\cap H^1_0(0,L)$ with homogeneous initial conditions and $t>0$.\\
	Let us introduce the Green operator $\mathcal{G}: H^{-1}(0,L)\rightarrow H^1_0(0,L)$ relative to $-\frac{\partial^2}{\partial x^2}$; then we have $\langle u,v\rangle_1=(\mathcal{G}^{1/2}u,\mathcal{G}^{1/2}v)_2$ for all $u, v\in H^{-1}(0,L)$.\\We omit at the moment the spatial dependence of the solutions; testing the two equations in (\ref{eq weak sist2sol}) respectively by $\varphi=\mathcal{G}\dot{w}$ and $\psi=\mathcal{G}\dot{\theta}$ we obtain
	\begin{equation}
	\begin{cases}
	&\dfrac{M}{2}\dfrac{d}{dt}||\dot{w}||^2_{H^{-1}}+\dfrac{EI}{2}\dfrac{d}{dt}||w||^2_{H^1}=\big(h_\alpha(w_1,\theta_1)-h_\alpha(w_2,\theta_{2}),\mathcal{G}^{1/2}\dot{w}\big)_2\\&\vspace{4mm}\hspace{0mm}+\big(h_\beta(w_1,\theta_1)-h_\beta(w_2,\theta_{2}),\mathcal{G}^{1/2}\dot{w}\big)_2\\
	&\dfrac{M\ell^2}{6}\dfrac{d}{dt}||\dot{\theta}||^2_{H^{-1}}+\dfrac{EJ}{2}\dfrac{d}{dt}||\theta||^2_{H^1}+\dfrac{GK}{2}\dfrac{d}{dt}||\theta||^2_{2}=\ell\big(h_\alpha(w_1,\theta_1),(\mathcal{G}\dot{\theta}\cos\theta_1)_x\big)_2+\\&-\ell\big(h_\alpha(w_2,\theta_2),(\mathcal{G}\dot{\theta}\cos\theta_2)_x\big)_2-\ell\big(h_\beta(w_1,\theta_1),(\mathcal{G}\dot{\theta}\cos\theta_1)_x\big)_2+\ell\big(h_\beta(w_2,\theta_2),(\mathcal{G}\dot{\theta}\cos\theta_2)_x\big)_2
	\end{cases}
	\label{eq sist sub e integr2}
	\end{equation}
	Now our aim is to find an upper bound for the right hand sides terms of (\ref{eq sist sub e integr2}). \\We observe that the nonlinearities $h_\alpha$, $h_\beta$, as defined in (\ref{g1,g2}), depend only on functions globally Lipschitzian. Indeed, introducing the integrand $\gamma(u):=\sqrt{1+[(u+y)_x]^2}$ of $\Gamma(u)$ and considering $\boldsymbol{\chi}(u)$, respectively as in (\ref{gamma}) and (\ref{chi}), we have that $\forall (x,t)\in(0,L)\times(0,\infty)\hspace{2mm} \exists\varrho:=\varrho(x,t)\in\big((w_1+\ell\sin\theta_1)_x,(w_2+\ell\sin\theta_2)_x\big)$ such that
	\small
	\begin{equation*}
	\begin{split}
	&|\gamma(w_1+\ell\sin\theta_{1})-\gamma(w_2+\ell\sin\theta_{2})|=\bigg|\frac{(\varrho+y_x)\big(w_1-w_2+\ell(\sin\theta_1-\sin\theta_2)\big)_x}{\sqrt{1+[\varrho+y_x]^2}}\bigg|\leq|w_x|+\ell|\theta_x|+\ell|\theta_{2x}\theta|,\\
	&\big|\boldsymbol{\chi}(w_1+\ell\sin\theta_1)-\boldsymbol{\chi}(w_2+\ell\sin\theta_2)\big|=\frac{\big|\big(w_1-w_2+\ell(\sin\theta_1-\sin\theta_2)\big)_x\big|}{(1+[\varrho+y_{x}]^2)^{\frac{3}{2}}}\leq|w_x|+\ell|\theta_x|+\ell|\theta_{2x}\theta|.
	\end{split}
	\end{equation*}
	\normalsize
	Then, recalling $h_\alpha$ as in (\ref{g1,g2}), we obtain
	\begin{equation}
	\begin{split}
	&|h_\alpha(w_1,\theta_1)-h_\alpha(w_2,\theta_{2})|\leq\big| H\xi\big\{\boldsymbol{\chi}(w_1+\ell\sin\theta_1)-\boldsymbol{\chi}(w_2+\ell\sin\theta_2)\big\}+\\&+\dfrac{AE_c}{L_c}\big\{\Gamma(w_1+\ell\sin\theta_1)\boldsymbol{\chi}(w_1+\ell\sin\theta_1)-\Gamma(w_2+\ell\sin\theta_2)\boldsymbol{\chi}(w_2+\ell\sin\theta_2)\big\}\big|\leq\\&\leq H\xi(|w_x|+\ell|\theta_x|+\ell|\theta_{2x}\theta|)+\dfrac{AE_c}{L_c}\big|\Gamma(w_1+\ell\sin\theta_1)\big[\boldsymbol{\chi}(w_1+\ell\sin\theta_1)-\boldsymbol{\chi}(w_2+\ell\sin\theta_2)\big]\\&+\boldsymbol{\chi}(w_2+\ell\sin\theta_2)\big[\Gamma(w_1+\ell\sin\theta_1)-\Gamma(w_2+\ell\sin\theta_2)\big]\big|\leq\\&\leq \bigg(H\overline{\xi}+\dfrac{AE_c}{L_c}\overline{C}\bigg)(|w_x|+\ell|\theta_x|+\ell|\theta_{2x}\theta|)+\dfrac{AE_c}{L_c}\int_{0}^{L}(|w_x|+\ell|\theta_x|+\ell|\theta_{2x}\theta|)dx,
	\label{stima h1}
	\end{split}
	\end{equation}
	in which $\overline{\xi}>\xi(x)$, see (\ref{xi}), and we have used again that $|\boldsymbol{\chi}|<1$ and $|\Gamma(u)|=\overline{C}$. Now considering  (\ref{stima h1}), applying the Schwartz and Young inequalities, it is possible to estimate the right hand side term of the first equation in (\ref{eq sist sub e integr2})
	\begin{equation}
	\big|\big(h_\alpha(w_1,\theta_1)-h_\alpha(w_2,\theta_{2}),\mathcal{G}^{1/2}\dot{w}\big)_2\big|\leq K_1\big(||\dot{w}||^2_{H^{-1}}+||w||^2_{H^1}+||\dot{\theta}||^2_{H^{-1}}+||\theta||^2_{H^1}+||\theta||^2_2\big).
	\label{stimaw}
	\end{equation}
	To obtain a similar result for the right hand side term of the second equation in (\ref{eq sist sub e integr2}) we need the following inequality 
	\begin{equation*}
	\begin{split}
	&\int_{0}^{L}|[\mathcal{G}\dot{\theta}(\cos\theta_1-\cos\theta_{2})]_x|dx\leq	\int_{0}^{L}\big(|\mathcal{G}^{1/2}\dot{\theta}(\cos\theta_1-\cos\theta_2)|+|\mathcal{G}\dot{\theta}\hspace{1mm}(\theta_{1x}\sin\theta_1-\theta_{2x}\sin\theta_2)|\big)dx\leq\\&\leq||\dot{\theta}||_{H^{-1}}||\cos\theta_1-\cos\theta_2||_2+||\dot{\theta}||_{H^*}||\theta||_{H^1}+||\dot{\theta}||_{H^*}||\theta_{2x}(\sin\theta_1-\sin\theta_2)||_{2}\leq\\&\leq K_2||\dot{\theta}||_{H^{-1}}\big(||\theta||_{H^1}+||\theta||_2\big),
	\end{split}
	\end{equation*}
	derived thanks to the Schwartz inequality, the embedding $H^{-1}\subset H^*$ and the Lipschitz property of the sine and cosine functions. Then the terms in the second equation of (\ref{eq sist sub e integr2}) are bounded
	\begin{equation}
	\begin{split}
	&\big|\big(h_\alpha(w_1,\theta_1),(\mathcal{G}\dot{\theta}\cos\theta_1)_x\big)_2-\big(h_\alpha(w_2,\theta_2),(\mathcal{G}\dot{\theta}\cos\theta_2)_x\big)_2\big|=\\&=\big|\big(h_\alpha(w_1,\theta_1)-h_\alpha(w_2,\theta_{2}),(\mathcal{G}\dot{\theta}\cos\theta_1)_x\big)_2+\big(h_\alpha(w_2,\theta_{2}),[\mathcal{G}\dot{\theta}(\cos\theta_1-\cos\theta_{2})]_x\big)_2\big|\leq\\&\leq K_3\big(||\dot{w}||^2_{H^{-1}}+||w||^2_{H^1}+||\dot{\theta}||^2_{H^{-1}}+||\theta||^2_{H^1}+||\theta||^2_2\big).
	\label{stimateta}
	\end{split}
	\end{equation}
	Next, integrating (\ref{eq sist sub e integr2}) in $s$ on $(0,t)$, adding the two left hand sides terms, we obtain a constant $\eta>0$ such that
	\begin{equation*}
	\begin{split}
	M||\dot{w}||^2_{H^{-1}}+EI||w||^2_{H^1}+\frac{M\ell^2}{3}||\dot{\theta}||^2_{H^{-1}}+EJ||\theta||^2_{H^1}+GK||\theta||^2_{2}\geq\\\hspace{20mm}\geq \eta\big(||\dot{w}||^2_{H^{-1}}+||w||^2_{H^1}+||\dot{\theta}||^2_{H^{-1}}+||\theta||^2_{H^1}+||\theta||^2_{2}\big).
	\end{split}
	\label{lower bound}
	\end{equation*}
	Hence, from (\ref{stimaw})-(\ref{stimateta}) and similar bounds for the terms in (\ref{eq sist sub e integr2}) involving the function $h_\beta$, we obtain $C>0$ such that
	\begin{equation*}
	\begin{split}
	&||\dot{w}||^2_{H^{-1}}+||w||^2_{H^1}+||\dot{\theta}||^2_{H^{-1}}+||\theta||^2_{H^1}+||\theta||^2_{2}\leq\\&\hspace{30mm}\leq C\int_{0}^{t}\big(||\dot{w}||^2_{H^{-1}}+||w||^2_{H^1}+||\dot{\theta}||^2_{H^{-1}}+||\theta||^2_{H^1}+||\theta||^2_{2}\big)ds.
	\end{split}
	\label{stima3}
	\end{equation*}
	Thanks to the Gronwall Lemma we have 
	\begin{equation*}
	\begin{split}
	&||\dot{w}(t)||^2_{H^{-1}}+||w(t)||^2_{H^1}+||\dot{\theta}(t)||^2_{H^{-1}}+||\theta(t)||^2_{H^1}+||\theta(t)||^2_{2}\leq\\&\hspace{27mm}\leq\big(||\dot{w}(0)||^2_{H^{-1}}+||w(0)||^2_{H^1}+||\dot{\theta}(0)||^2_{H^{-1}}+||\theta(0)||^2_{H^1}+||\theta(0)||^2_{2}\big)e^{Ct},
	\end{split}
	\label{gr}
	\end{equation*}
	and this fact ensures
	\begin{equation*}
	\begin{split}
	&||\dot{w}||^2_{H^{-1}}+||w||^2_{H^1}+||\dot{\theta}||^2_{H^{-1}}+||\theta||^2_{H^1}+||\theta||^2_{2}=0 \hspace{5mm}\forall t\in [0,T].
	\end{split}
	\label{fine}
	\end{equation*}
	In this way the uniqueness of the weak solution $(w,\theta)\in Z^2_T$, where $ Z_T:=C^0\big([0,T];H^1_0(0,L)\big)\cap C^1\big([0,T];H^{-1}(0,L)\big)$, is obtained. Thanks to the regularity of $(w,\theta)$ and the fact that $X_T\subset Z_T$ we have a unique weak solution $(w,\theta)\in X^2_T$, satisfying the initial conditions (\ref{ics2}). This fact completes the proof of Theorem \ref{teo1}.
	\begin{flushright}
		$\square$
	\end{flushright}
	\section{Conclusions}
	In this paper we have presented an isolated model for suspension bridges with deformable cables and rigid hangers, inspired by the Melan equation. Differently from it, the system of partial differential equations is derived from variational principles; we have considered two degrees of freedom, the vertical displacement and the torsional rotation of the deck, to see when the torsional instability phenomena arise.\\
	The system obtained is nonlinear due to the geometric configuration of the cables and to the rotation of the deck, in particular we avoid assumptions on small rotations of the deck; moreover, we involve in the model not only the torsional effects on the deck due to de Saint Venant theory, but also those more precise introduced by Vlasov theory.\\
	The presence of two cables linked to a single deck produces a problem more complex than the original Melan equation (single cable-beam system) and gives two strongly coupled equations of the motion.\\
	Adopting the Galerkin procedure we proved the existence of a solution, while we proved its uniqueness testing the equations with the Green function applied to the time derivative of the solutions; we proposed the complete proof since it is non-standard due to the presence of these nonlinearities.\\
	We have shown some numerical experiments on this system, considering 10 longitudinal modes interacting with 4 torsional modes. The results show that, exciting distinct longitudinal modes, there are different thresholds of torsional instability; this fact reveals that the origin of the instability is structural, as \cite{berchio, Gazzola Torsion, falocchi}. It is clear that in absence of wind, the deck does not move, but when the wind hits a bluff body vertical oscillations begin, due to the vortex effect. When we apply the initial condition $\overline{w}_9^0 =3.87$m a periodic motion rises on this longitudinal component and it is maintained in amplitude by a somehow perfect equilibrium between the input of energy and the structural damping. This is why, at least as a first step, it appears reasonable to consider isolated systems.\\
	In the paper an analysis of sensitivity with respect to the mechanical parameters involved in the system is performed; we have considered the parameters of the TNB as "basic situation" and then, modifying them, we have analyzed how the response of the system changes in terms of torsional instability. These observations can give some hints to bridges designer.\\
	A possible development in this field is the study of the model involving the slackening mechanisms of the hangers. A natural way to introduce this assumption is to set up a model in which the degrees of freedom pass from 2 to 4; this choice implies a further effort in terms of modeling and computation. For these reasons, in a future work, we are planning to set up a model with 2 degrees of freedom and, acting on the cable shape, we will model indirectly the slackening of the hangers.

\end{document}